\documentclass[a4paper,11pt]{amsart}

\usepackage{graphicx}
\usepackage{mathptmx}
\usepackage{amsmath}
\usepackage{amssymb}
\usepackage{enumitem}
\usepackage{xcolor}

\newmuskip\pFqmuskip

\newcommand*\pFq[6][8]{%
  \begingroup 
  \pFqmuskip=#1mu\relax
  \mathcode`=\string"8000
  \begingroup\lccode`\~=`\,
  \lowercase{\endgroup\let~}\pFqcomma
  F^{#2}_{#3}{\left(\genfrac..{0pt}{}{#4}{#5}\bigg|#6\right)}%
  \endgroup
}
\newcommand{\pFqcomma}{\mskip\pFqmuskip}

\newtheorem{theorem}{Theorem}[section]
\newtheorem{lemma}[theorem]{Lemma}
\newtheorem{corollary}[theorem]{Corollary}

\newtheorem{remark}[theorem]{Remark}

\begin{document}

\title[]{Combinatorial identities related to degenerate Stirling numbers of the second kind}

\author{Dae San  Kim}
\address{Department of Mathematics, Sogang University, Seoul 121-742, Republic of Korea}
\email{dskim@sogang.ac.kr}
\author{Taekyun  Kim}
\address{Department of Mathematics, Kwangwoon University, Seoul 139-701, Republic of Korea}
\email{tkkim@kw.ac.kr}

\subjclass[2010]{11B68; 11B73; 11B83}
\keywords{degenerate Stirling numbers of the second kind; degenerate Bernoulli numbers of order $\alpha$}

\begin{abstract}
The study of degenerate versions of certain special polynomials and numbers, which was initiated by Carlitz's work on degenerate Euler and degenerate Bernoulli polynomials, has recently seen renewed interest among mathematicians.
The aim of this paper is to study some properties, certain identities, recurrence relations and explicit expressions for degenerate Stirling  numbers of the second kind, which are a degenerate version of the Stirling numbers of the second kind. These numbers appear very frequently when we study various degenerate versions of many special polynomials and numbers. Especially, we consider some closely related polynomials and power series in connection with a degenerate version of Euler's formula for the Stirling numbers of the second kind.
\end{abstract}

\maketitle

\markboth{\centerline{\scriptsize Combinatorial identities related to degenerate Stirling numbers of the second kind}}
{\centerline{\scriptsize Dae San Kim and Taekyun Kim}}

\section{Introduction} 
The Stirling number of the second ${n \brace k}$ is the number of ways to partition a set of $n$ objects into $k$ nonempty subsets. Stirling numbers of the second kind appear as the coefficients in the expansion of the sequence  of polynomials ${\bf{P}}=\{x^n\}_{n=0}^{\infty}$ in terms of falling factorials (see \eqref{1}, \eqref{1-1}, \eqref{7-1}). This is generalized to any sequence of polynomials ${\bf{P}}=\{p_{n}(x)\}_{n=0}^{\infty}$ satisfying deg\,$p_{n}(x)=n,\, p_{0}(x)=1$ (see [9]). Indeed, for any such a sequence of polynomials ${\bf{P}}=\{p_{n}(x)\}_{n=0}^{\infty}$, the Stirling numbers of the second kind associated with $\bf{P}$, $S_2(n,k;{\bf{P}})$, are defined as the coefficients when we expand $p_{n}(x)$ in terms of $(x)_{k}$:
\begin{equation*}
p_{n}(x)=\sum_{k=0}^{n}S_2(n,k;{\bf{P}})(x)_{k}.
\end{equation*}
The degenerate Stirling number ${n \brace k}_{\lambda}$ (see \eqref{4}), the (unsigned) Lah number $L(n,k)$ and Gould-Hopper number $G(n,k;r,s)$ arise as the coefficients respectively in the expansion of the degenerate falling factorial $(x)_{n,\lambda}$, the rising factorial $\langle{x\rangle}_{n}$ and $(rx+s)_{n}$ in terms of $(x)_{k}$, where $r \ne 0$, and $\langle{x\rangle}_{0}=1,
\langle{x\rangle}_{n}=x(x+1)\cdots(x+n-1)$, for $n \ge 1$. Furthermore, the Stirling numbers of the first kind associated with $\bf{P}$, $S_{1}(n,k;{\bf{P}})$, are defined as $(x)_{n}=\sum_{k=0}^{n}S_{1}(n,k;{\bf{P}})p_{k}(x)$. The Stirling numbers of both kinds associated with any sequence of polynomials were investigated with the help of umbral calculus. They are illustrated for many examples of ${\bf{P}}=\{p_{n}(x)\}_{n=0}^{\infty}$, together with their orthogonality and inverse relations in each case (see [9]).\par
Another generalization of the Stirling numbers of both kinds are the recent study of the probabilistic Stirling numbers of both kinds. Assume that $Y$ is a random variable whose moment generating function $\mathbb{E}[e^{tY}]=\sum_{n=0}^{\infty}\mathbb{E}(Y^{n})\frac{t^n}{n!}$ exists in a neighborhood of the origin. Then the probabilistic Stirling numbers of the second kind associated with $Y$ are defined as (see [21])
\begin{equation}
\frac{1}{k!}(\mathbb{E}[e^{tY}]-1)^{k}=\sum_{n=k}^{\infty}S_{2}^{Y}(n,k)\frac{t^n}{n!}.
\end{equation}
We note that if $Y=1$ then $S_{2}^{1}(n,k)={n \brace k}$. Furthermore, the probabilistic degenerate Stirling numbers of the second kind associted with $Y$ are defined as $\frac{1}{k!}(\mathbb{E}[e_{\lambda}^{Y}(t)]-1)^{k}=\sum_{n=k}^{\infty}S_{2, \lambda}^{Y}(n,k)\frac{t^n}{n!}$. Here we see that $S_{2,\lambda}^{1}(n,k)={n \brace k}_{\lambda}$, for $Y=1$. The reader may refer to the papers [2,12,13,15,16,21,26 and the references therein] for the recent studies on probabilistic extensions of certain special polynomials and numbers.\par
It is Carlitz [5,6] who initiated a study of degenerate versions of the ordinary Bernoulli and Euler polynomials and numbers, namely the degenerate Bernoulli and Euler polynomials and numbers. The study of degenerate versions of certain special polynomials and numbers has recently seen renewed interest among mathematicians. They have been explored by various methods, which include generating functions, umbral calculus techniques, combinatorial methods, $p$-adic analysis, differential equations, special functions, probability theory and analytic number theory. We let the reader refer to the papers [5,10,11,14,16-18,20-25, 29,32 and the references therein] for recent results on degenerate versions of certain special numbers and polynomials. It is remarkable that studying degenerate versions is not only limited to polynomials but also extended to transcendental functions. Indeed, the degenerate gamma functions were introduced in connection with degenerate Laplace transforms (see [14]). It is also noteworthy that the degenerate umbral calculus is introduced as a degenerate version of the classical umbral calculus (see [10]). \par
The degenerate Stirling numbers of the second kind ${n \brace k}_{\lambda}$ (see \eqref{4}, \eqref{5}, \eqref{6}, \eqref{7}) are a degenerate version of the ordinary Stirling numbers of the second kind and appear very frequently when we study various degenerate versions of many special numbers and polynomials. The aim of this paper is to study some properties, certain identities, recurrence relations and explicit expressions related to the degenerate Stirling numbers of the second kind. Especially, we study the closely related polynomials and power series, $S_{n}(x, r|\lambda), K_{r}(x|\lambda), S_{n,r}(x|\lambda)$, and $T_{n}(x, r|\lambda)$, in connection with the degenerate version of Euler's formula for the Stirling numbers of the second kind in \eqref{6}. The reader may refer to [2,4,5,8,9,11,17-19,21,22 and the references therein] for Stirling numbers and [1,3,7,27,28,30,31] as general references. In the rest of this section, we recall the facts that are needed throughout this paper. \par

\vspace{0.1in}

For $n\ge 0$, the Stirling numbers of the second kind are defined by 
\begin{equation}
x^{n}=\sum_{k=0}^{n}{n \brace k}(x)_{k},\quad (\mathrm{see}\ [1,2,6,27]), \label{1}	
\end{equation}
where the falling factorials $(x)_{n}$ are given by
\begin{equation}
(x)_{0}=1,\ (x)_{n}=x(x-1)(x-2)\cdots(x-n+1),\ (n\ge 1). \label{1-1}
\end{equation} 
For any nonzero $\lambda\in\mathbb{R}$, the degenerate exponentials are defined by 
\begin{equation}
e_{\lambda}^{x}(t)=\sum_{k=0}^{\infty}(x)_{k,\lambda}\frac{t^{k}}{k!}=(1+\lambda t)^{\frac{x}{\lambda}},\quad  e_{\lambda}(t)=e_{\lambda}^{1}(t),\quad (\mathrm{see}\ [9-22]), \label{2}	
\end{equation}
where the degenerate falling factorials $(x)_{n,\lambda}$ are given by
\begin{equation}
(x)_{0,\lambda}=1,\ (x)_{n,\lambda}=x(x-\lambda)(x-2\lambda)\cdots\big(x-(n-1)\lambda\big),\quad (n\ge 1). \label{3}
\end{equation}
Note that 
\begin{displaymath}
\lim_{\lambda\rightarrow 0}e_{\lambda}^{x}(t)=e^{xt}. 
\end{displaymath} 
Recently, the degenerate Stirling numbers of the second kind are defined by 
\begin{equation}
(x)_{n,\lambda}=\sum_{k=0}^{n}{n \brace k}_{\lambda}(x)_{k},\quad (n\ge 0),\quad (\mathrm{see}\ [9,15,17]).\label{4}
\end{equation}
From \eqref{4}, we can show 
\begin{equation}
\frac{1}{k!}\big(e_{\lambda}(t)-1\big)^{k}=\sum_{n=k}^{\infty}{n \brace k}_{\lambda}\frac{t^{n}}{n!},\quad (k\ge 0),\quad (\mathrm{see}\ [9,12,18]), \label{5}
\end{equation}
and the degenerate version of Euler's formula given by
\begin{equation}
\frac{1}{k!}\sum_{j=0}^{k}\binom{k}{j}(-1)^{k-j}(j)_{n,\lambda}={n \brace k}_{\lambda},\quad (n,k \ge 0),\label{6}
\end{equation}
with the understanding that ${n \brace k}_{\lambda}=0$, for $k > n$. \\
By \eqref{1} and \eqref{4}, we have $\lim_{\lambda\rightarrow 0}{n \brace k}_{\lambda}={n \brace k},\ (n\ge k\ge 0).$ It is not difficult to show that 
\begin{equation}
{n+1 \brace k}_{\lambda}={n \brace k-1}_{\lambda}+(k-n\lambda){n \brace k}_{\lambda},\quad (n\ge k\ge 1). \label{7}
\end{equation} \\
By taking $\lambda \rightarrow 0$ in \eqref{5}, \eqref{6} and \eqref{7}, the Stirling numbers of the second kind ${n \brace k}$ are alternatively given by 
\begin{align}
&\frac{1}{k!}\big(e^{t}-1\big)^{k}=\sum_{n=k}^{\infty}{n \brace k}\frac{t^{n}}{n!},\quad (k\ge 0), \nonumber\\
&\frac{1}{k!}\sum_{j=0}^{k}\binom{k}{j}(-1)^{k-j}j^{n}={n \brace k},\quad (n,k \ge 0), \label{7-1}\\
&{n+1 \brace k}={n \brace k-1}+k{n \brace k},\quad (n\ge k\ge 1).\nonumber
\end{align} \par
In [5], Carlitz introduced the degenerate Bernoulli numbers of order $\alpha$ which are given by 
\begin{equation}
\bigg(\frac{t}{e_{\lambda}(t)-1}\bigg)^{\alpha}=\sum_{n=0}^{\infty}\beta_{n,\lambda}^{(\alpha)}\frac{t^{n}}{n!}.\label{8}
\end{equation}
When $\alpha=1, \beta_{n,\lambda}=\beta_{n,\lambda}^{(1)}$ are called the degenerate Bernoulli numbers. \\
Note that $\displaystyle\lim_{\lambda\rightarrow 0}\beta_{n,\lambda}^{(\alpha)}=B_{n}^{(\alpha)},\ (n\ge 0)\displaystyle$, where $B_{n}^{(\alpha)}$ are the Bernoulli numbers of order $\alpha$, given by 
\begin{equation}
\bigg(\frac{t}{e^{t}-1}\bigg)^{\alpha}=\sum_{n=0}^{\infty}B_{n}^{(\alpha)}\frac{t^{n}}{n!},\quad (\mathrm{see}\ [1-31]). \label{9}
\end{equation} \par
As usual, the general binomial coefficients $\binom{x}{n}$ are defined as follows:
\begin{equation}
\binom{x}{n}=\frac{(x)_{n}}{n!},\quad(n \ge 0). \label{9-1}
\end{equation}
The facts will be repeatedly used throughout this paper.
\begin{lemma}
For any nonnegative integer $n$, the following hold true.
\begin{flalign*}
&(a)\,\, (-x)_{n,\lambda}=(-1)^{n}(x)_{n,-\lambda},\\
&(b)\,\, \binom{x}{n}=(-1)^{n}\binom{n-x-1}{n}, \\
&(c)\,\, \binom{x+y}{n}=\sum_{k=0}^{n}\binom{x}{k}\binom{y}{n-k},\\
&(d)\,\,(x+y)_{n,\lambda}=\sum_{k=0}^{n}\binom{n}{k}(x)_{k,\lambda}(y)_{n-k,\lambda},\\
&(e)\,\,\sum_{k=0}^{n}(-1)^{k}\binom{x}{k}=(-1)^{n}\binom{x-1}{n}. &&
\end{flalign*}
\begin{proof}
(a) and (b) follow from definitions (see \eqref{1-1}, \eqref{3}, \eqref{9-1}). (c) is the Vandermonde convolution formula which follows from $(1+t)^{x+y}=(1+t)^{x}(1+t)^{y}$. Similarly to (c), (d) follows from $e_{\lambda}^{x+y}(t)=e_{\lambda}^{x}(t)e_{\lambda}^{y}(t)$ (see \eqref{2}). The left hand of (e) becomes a telescoping sum upon replacing $\binom{x}{k}$ by $\binom{x-1}{k-1}+\binom{x-1}{k}=\binom{x}{k}$.
\end{proof}
\end{lemma}
We use the notations $D_{x}=\frac{d}{dx},\, D_{z}=\frac{d}{dz},\, D_{x}^{n}=(\frac{d}{dx})^{n}$,\,  and  $D_{z}^{n}=(\frac{d}{dz})^{n}$. Assume that both $f(z)$ and $z=\phi(x)$ are real or complexed valued functions which are $n$-times differentiable. Then the Hoppe's generalized chain rule says that (see [27], Theorem 8.1, p.105) 
\begin{equation}
D_{x}^{n}f(z)=\sum_{k=0}^{n}D_{z}^{k}f(z)\frac{(-1)^k}{k!}\sum_{j=0}^{k}(-1)^{j}\binom{k}{j}z^{k-j}D_{x}^{n}z^{j}. \label{9-2}
\end{equation}
As it is mentioned in ([27], (8.24), p.108), by applying \eqref{9-2} to $f(z)=z^{-\alpha}$ we get the following lemma. 
\begin{lemma}
Assume that $z$ is a function of $x$ that is $n$-times differentiable and $\alpha$ is a real or complex number. Then we have
\begin{equation*}
D_{x}^{n}z^{-\alpha}=\alpha \binom{\alpha+n}{n}\sum_{j=0}^{n}(-1)^{j}\binom{n}{j}\frac{z^{-\alpha-j}}{\alpha+j}D_{x}^{n}z^{j}.
\end{equation*}
\end{lemma}

\section{Combinatorial identities related to degenerate Stirling numbers of the second kind} 
For $n,\alpha\in\mathbb{N}$, by \eqref{7}, we get 
\begin{equation}
\begin{aligned}
	&\alpha!{n \brace \alpha}_{\lambda}=\alpha!\big(\alpha-(n-1)\lambda\big){n-1 \brace \alpha}_{\lambda}+\alpha!{n-1 \brace \alpha-1}_{\lambda} \\
	&=\alpha!\alpha{n-1 \brace \alpha}_{\lambda}+\alpha (\alpha-1)!{n-1 \brace \alpha-1}_{\lambda}-\alpha(\alpha-1)!(n-1)\lambda {n-1 \brace \alpha}_{\lambda}. 
\end{aligned}	\label{10}
\end{equation}
Replacing $\alpha$ by $\alpha-1$ and $n$ by $n-1$ in \eqref{10}, we get 
\begin{equation}
\begin{aligned}
	&(\alpha-1)!{n-1 \brace \alpha-1}_{\lambda}=(\alpha-1)(\alpha-1)!{n-2 \brace \alpha-1}_{\lambda}\\
	&+(\alpha-1)(\alpha-2)!{n-2 \brace \alpha-2}_{\lambda} 
-(\alpha-1)(\alpha-2)!(n-2)\lambda{n-2 \brace \alpha-1}_{\lambda}. \label{11}
\end{aligned}	
\end{equation}
By \eqref{10} and \eqref{11}, we get 
\begin{align}
&\alpha!{n\brace \alpha}_{\lambda}=\alpha!\alpha{n-1\brace \alpha}_{\lambda}-\alpha(\alpha-1)!(n-1)\lambda	{n-1 \brace \alpha}_{\lambda}\label{12}\\
&+\alpha(\alpha-1)(\alpha-1)!{n-2 \brace \alpha-1}_{\lambda}+\alpha(\alpha-1)(\alpha-2)!{n-2 \brace \alpha-2}_{\lambda}\nonumber \\
&-\alpha(\alpha-1)(\alpha-2)!(n-2)\lambda{n-2 \brace \alpha-1}_{\lambda}.\nonumber
\end{align}
Replacing $n$ by $n-2$ and $\alpha$ by $\alpha-2$ in \eqref{10}, we have 
\begin{align}
&(\alpha-3)!{n-2 \brace \alpha-2}=(\alpha-2)(\alpha-2)!{n-3 \brace \alpha-2}_{\lambda}\label{13}\\
&+(\alpha-2)(\alpha-3)!{n-3 \brace \alpha-3}_{\lambda}-(\alpha-2)(\alpha-3)!(n-3)\lambda{n-3 \brace \alpha-2}_{\lambda}. \nonumber
\end{align}
From \eqref{12} and \eqref{13}, we note that 
\begin{align}
&\alpha!{n \brace \alpha}_{\lambda}=\alpha!\alpha{n-1 \brace \alpha}_{\lambda}-\alpha (\alpha-1)!(n-1)\lambda	{n-1 \brace \alpha}_{\lambda}\label{14}\\
&+\alpha(\alpha-1)(\alpha-1)!{n-2 \brace \alpha-1}_{\lambda}+\alpha(\alpha-1)(\alpha-2)(\alpha-2)!{n-3 \brace \alpha-2}_{\lambda}\nonumber \\
&+\alpha(\alpha-1)(\alpha-2)(\alpha-3)!{n-3 \brace \alpha-3}_{\lambda}-\alpha(\alpha-1)(\alpha-2)!(n-2)\lambda{n-2 \brace \alpha-1}_{\lambda}\nonumber\\
&-\alpha(\alpha-1)(\alpha-2)(\alpha-3)!(n-3)\lambda{n-3 \brace \alpha-2}_{\lambda}.\nonumber
\end{align}
Continuing this process, we have 
\begin{align}
\alpha!{n \brace \alpha}_{\lambda}&=\sum_{k=0}^{r}\binom{\alpha}{k+1}(k+1)!(\alpha-k)!{n-1-k \brace \alpha-k}_{\lambda} \label{15}	\\
&\quad +\binom{\alpha}{r+1}(r+1)!(\alpha-r-1)!{n-r-1 \brace \alpha-r-1}_{\lambda}\nonumber \\
&-\lambda\sum_{k=0}^{r}\binom{\alpha}{k+1}(k+1)!(n-k-1){n-1-k \brace \alpha-k}_{\lambda}(\alpha-k-1)!\nonumber \\
&=\sum_{k=0}^{r}\binom{\alpha}{k+1}(k+1)!(\alpha-k-1)!\big(\alpha-k-\lambda(n-k-1)\big){n-1-k \brace \alpha-k}_{\lambda}\nonumber\\
&+\binom{\alpha}{r+1}(r+1)!(\alpha-r-1)!{n-r-1 \brace \alpha-r-1}_{\lambda}. \nonumber
\end{align}
For $n-\alpha>0$, we have ${n-\alpha \brace 0}_{\lambda}=0$. Let $n > \alpha$, and let $r=\alpha-1$. Then, by \eqref{15}, we get 
\begin{align}
\alpha!{n \brace \alpha}_{\lambda}&=\sum_{k=0}^{\alpha-1}\binom{\alpha}{k+1}(k+1)!(\alpha-k-1)!\big(\alpha-k-\lambda(n-k-1)\big){n-1-k \brace \alpha-k}_{\lambda}\label{16}\\
&\quad +\binom{\alpha}{\alpha}\alpha!0!{n-\alpha \brace 0}_{\lambda}\nonumber \\
&=\sum_{k=1}^{\alpha}\binom{\alpha}{k}k!(\alpha-k)!\big(\alpha-k+1-\lambda(n-k)\big){n-k \brace \alpha-k+1}_{\lambda}\nonumber\\
&=\alpha!\sum_{k=1}^{\alpha}\big(\alpha-k+1-\lambda(n-k)\big){n-k \brace \alpha-k+1}_{\lambda}.\nonumber
\end{align}
Therefore, by \eqref{16}, we obtain the following theorem. 
\begin{theorem}
For $\alpha,n\in\mathbb{N}$ with $n>\alpha$, we have 
\begin{displaymath}
{n \brace \alpha}_{\lambda}=\sum_{k=1}^{\alpha}\big(\alpha-k+1-\lambda(n-k)\big){n-k \brace \alpha-k+1}_{\lambda}.
\end{displaymath}
\end{theorem}
\begin{remark} 
\begin{align*}
	{n \brace n-1}_{\lambda}&=\sum_{k=1}^{n-1}\big(n-k-\lambda(n-k)\big){n-k \brace n-k}_{\lambda} \\
	&=(1-\lambda)\sum_{k=1}^{n-1}(n-k)=(1-\lambda)\binom{n}{2}, 
\end{align*}
and 
\begin{align*}
	{n \brace n-2}_{\lambda}&=\sum_{k=1}^{n-2}\big((n-k-1)-\lambda(n-k)\big){n-k \brace n-k-1}_{\lambda}\\
	&=(1-\lambda)\sum_{k=1}^{n-2}\big((n-k-1)-\lambda(n-k)\big)\binom{n-k}{2}\nonumber \\
	&=(1-\lambda)\sum_{k=1}^{n-2}(1-\lambda)(n-k-1)\binom{n-k}{2}-\lambda(1-\lambda)\sum_{k=1}^{n-2}\binom{n-k}{2}\\
	&=(1-\lambda)^{2}\sum_{k=1}^{n-2}\frac{k^{2}(k+1)}{2}-\lambda(1-\lambda)\sum_{k=1}^{n-2}\frac{k(k+1)}{2}\\
	&=\frac{(1-\lambda)^{2}}{2}\bigg(\sum_{k=1}^{n-2}k^{3}+\sum_{k=1}^{n-2}k^{2}\bigg)-\frac{\lambda(1-\lambda)}{2}\bigg(\sum_{k=1}^{n-2}k^{2}+\sum_{k=1}^{n-2}k\bigg)\\
	&=\frac{1}{24}(1-\lambda)(3n-5-\lambda(3n-1))n(n-1)(n-2). 
\end{align*}
\end{remark}
For $n\ge 0$, by \eqref{4} and Lemma 1.1 (a), (b), (c), we get 
\begin{align}
&(-1)^{n}\sum_{\alpha=0}^{n}{n \brace \alpha}_{-\lambda}\binom{x}{\alpha}\alpha!=(-1)^{n}(x)_{n,-\lambda}=(-x)_{n,\lambda}=\sum_{j=0}^{n}\binom{-x}{j}j!{n \brace j}_{\lambda}\label{17} \\
&=\sum_{j=0}^{n}(-1)^{j}\binom{x+j-1}{j}j!{n \brace j}_{\lambda}=\sum_{j=0}^{n}(-1)^{j}j!{n \brace j}_{\lambda}\sum_{\alpha=0}^{j}\binom{x}{\alpha}\binom{j-1}{j-\alpha}\nonumber\\
&=\sum_{\alpha=0}^{n}\binom{x}{\alpha}\sum_{j=\alpha}^{n}(-1)^{j}j!{n \brace j}_{\lambda}\binom{j-1}{\alpha-1}.\nonumber
\end{align}
Therefore, by comparing the coefficients on both sides of \eqref{17}, we obtain the following theorem. 
\begin{theorem}
Let $\alpha, n$ be nonnegative integers with $n\ge \alpha$. Then we have 
\begin{displaymath}
\alpha!{n \brace \alpha}_{-\lambda}=\sum_{j=\alpha}^{n}(-1)^{n-j}j!{n \brace j}_{\lambda}\binom{j-1}{\alpha-1}. 
\end{displaymath}
\end{theorem}
Now, we observe by (4) and Lemma 1.1 (c) that 
\begin{align}
(-1)^{n}(-x)_{n,\lambda}&=x(x+\lambda)\cdots\big(x+(n-1)\lambda\big)=\big(x+(n-1)\lambda\big)_{n,\lambda}\label{18} \\
&=\sum_{j=0}^{n}{n \brace j}_{\lambda}\binom{x+(n-1)\lambda}{j}j!\nonumber\\
&=\sum_{j=0}^{n}{n \brace j}_{\lambda}j!\sum_{\alpha=0}^{j}\binom{x}{\alpha}\binom{(n-1)\lambda}{j-\alpha}\nonumber\\
&=\sum_{\alpha=0}^{n}\sum_{j= \alpha}^{n}{n \brace j}_{\lambda}j!\binom{(n-1)\lambda}{j-\alpha}\binom{x}{\alpha}.\nonumber
\end{align}
Therefore, by \eqref{17} and \eqref{18}, we obtain the following theorem. 
\begin{theorem}
Let $\alpha,n$ be nonnegative integer with $n\ge \alpha$. Then we have 
\begin{displaymath}
\sum_{j=\alpha}^{n}{n \brace j}_{\lambda}j!\binom{(n-1)\lambda}{j-\alpha}=\alpha!{n \brace \alpha}_{-\lambda}=\sum_{j=\alpha}^{n}(-1)^{n-j}j!{n \brace j}_{\lambda}\binom{j-1}{\alpha-1}. 
\end{displaymath}
\end{theorem}
We define the difference operator $\triangle_{x,h}$ by 
\begin{equation}
\triangle_{x,h}f(x)=\frac{1}{h}\big(f(x+h)-f(x)\big).\label{19}
\end{equation}
Then, by \eqref{19}, we get 
\begin{equation}
\triangle_{x,h}^{n}f(x)=\frac{1}{h^{n}}\sum_{k=0}^{n}\binom{n}{k}(-1)^{n-k}f(hk+x),\quad (n\in\mathbb{N}).\label{20}
\end{equation}
Let $f(x)=(x)_{p,\lambda},\ (p\in\mathbb{N})$. From \eqref{6} and \eqref{20} and using Lemma 1.1 (d), we have 
\begin{align}
\triangle_{x,h}^{n}(x)_{p,\lambda}&=\frac{1}{h^{n}}\sum_{k=0}^{n}\binom{n}{k}(-1)^{n-k}(x+kh)_{p,\lambda}\label{21}\\
&=\frac{1}{h^{n}}\sum_{k=0}^{n}(-1)^{n-k}\binom{n}{k}\sum_{r=0}^{p}(kh)_{r,\lambda}(x)_{p-r,\lambda}\binom{p}{r}\nonumber\\
&=\frac{1}{h^{n}}\sum_{r=0}^{p}\binom{p}{r}(x)_{p-r,\lambda}h^{r}\sum_{k=0}^{n}\binom{n}{k}(-1)^{n-k}(k)_{r,\frac{\lambda}{h}}\nonumber\\
&=\sum_{r=0}^{p}\binom{p}{r}(x)_{p-r,\lambda}h^{r-n}n!{r \brace n}_{\frac{\lambda}{h}}. \nonumber
\end{align}
Therefore, by \eqref{21}, we obtain the following theorem. 
\begin{theorem}
For $n,p \ge 0$, we have 
\begin{displaymath}
\frac{1}{n!}\triangle_{x,h}^{n}(x)_{p,\lambda}=\sum_{r=0}^{p}\binom{p}{r}(x)_{p-r,\lambda}h^{r-n}{r \brace n}_{\frac{\lambda}{h}}. 
\end{displaymath}
\end{theorem}
For $n,j\ge 0$, by \eqref{6} and Lemma 1.1 (a), (d), we get 
\begin{align}
j!{n \brace j}_{\lambda}&=\sum_{k=0}^{j}(-1)^{k-j}\binom{j}{k}(k)_{n,\lambda}=\sum_{k=0}^{j}(-1)^{k}\binom{j}{k}(j-k)_{n,\lambda}\label{22} \\
&=\sum_{k=0}^{j}(-1)^{k}\binom{j}{k}\sum_{\alpha=0}^{n}\binom{n}{\alpha}(j)_{n-\alpha,\lambda}(-k)_{\alpha,\lambda}\nonumber\\
&=\sum_{k=0}^{j}(-1)^{k}\binom{j}{k}\sum_{\alpha=0}^{n}\binom{n}{\alpha}(j)_{n-\alpha,\lambda}(-1)^{\alpha}(k)_{\alpha,-\lambda}\nonumber\\
&=\sum_{\alpha=0}^{n}\binom{n}{\alpha}(-1)^{\alpha+j}(j)_{n-\alpha,\lambda}j!\frac{1}{j!}\sum_{k=0}^{j}\binom{j}{k}(-1)^{j-k}(k)_{\alpha,-\lambda}\nonumber\\
&=\sum_{\alpha=0}^{n}\binom{n}{\alpha}(-1)^{\alpha+j}(j)_{n-\alpha,\lambda}j!{\alpha \brace j}_{-\lambda}. \nonumber
\end{align} \\
Therefore, by \eqref{22}, we obtain the following theorem. 
\begin{theorem}
	For $n,j\ge 0$ with $n\ge j$, we have 
	\begin{displaymath}
		{n \brace j}_{\lambda}=\sum_{\alpha=0}^{n}\binom{n}{\alpha}(-1)^{\alpha+j}{\alpha \brace j}_{-\lambda}(j)_{n-\alpha,\lambda}.
	\end{displaymath}
\end{theorem}
For $\alpha,r\in\mathbb{N}$, by \eqref{5}, we get 
\begin{align}
\binom{r+\alpha}{r}\sum_{n=r+\alpha}^{\infty}{n \brace r+\alpha}_{\lambda}\frac{t^{n}}{n!}&=\frac{1}{r!}\big(e_{\lambda}(t)-1\big)^{r}\frac{1}{\alpha!}\big(e_{\lambda}(t)-1\big)^{\alpha}\label{23} \\
&=\sum_{k=r}^{\infty}{k \brace r}_{\lambda}\frac{t^{k}}{k!}\sum_{l=\alpha}^{\infty}{l \brace \alpha}_{\lambda}\frac{t^{l}}{l!}\nonumber \\
&=\sum_{n=\alpha+r}^{\infty}\sum_{k=r}^{n}\binom{n}{k}{k \brace r}_{\lambda}{n-k \brace \alpha}_{\lambda}\frac{t^{n}}{n!}.\nonumber
\end{align}
Thus, for $n \ge r+\alpha$, we have 
\begin{equation}
\binom{r+\alpha}{r}{n \brace r+\alpha}_{\lambda}=\sum_{k=r}^{n}\binom{n}{k}{k \brace r}_{\lambda}{n-k \brace \alpha}_{\lambda}.\label{24}	
\end{equation}
For $n,\alpha\in\mathbb{N}$, by \eqref{4} and Lemma 1.1 (d), we have 
\begin{align}
&\sum_{j=0}^{n}(-1)^{j}\binom{n}{j}(z-j)_{n+\alpha,\lambda}=\sum_{j=0}^{n}(-1)^{j}\binom{n}{j}(z-n+n-j)_{n+\alpha,\lambda}\label{25}\\
&=\sum_{j=0}^{n}\binom{n}{j}(-1)^{j}\sum_{k=0}^{n+\alpha}\binom{n+\alpha}{k}(n-j)_{n+\alpha-k,\lambda}(z-n)_{k,\lambda}\nonumber\\
&=\sum_{j=0}^{n}\binom{n}{j}(-1)^{j}\sum_{k=0}^{n+\alpha}\binom{n+\alpha}{k}(n-j)_{n+\alpha-k,\lambda}\sum_{r=0}^{k}{k \brace r}_{\lambda}\binom{z-n}{r}r!\nonumber\\
&=\sum_{r=0}^{n+\alpha}\binom{z-n}{r}\sum_{k=r}^{n+\alpha}\binom{n+\alpha}{k}r!{k \brace r}_{\lambda}\sum_{j=0}^{n}\binom{n}{j}(-1)^{j}(n-j)_{n+\alpha-k,\lambda}.\nonumber
\end{align}
We note from \eqref{6} that 
\begin{align}
\sum_{j=0}^{n}\binom{n}{j}(-1)^{j}(n-j)_{n+\alpha-k,\lambda}&=n! \frac{1}{n!}\sum_{j=0}^{n}\binom{n}{n-j}(-1)^{n-(n-j)}(n-j)_{n+\alpha-k,\lambda} \label{25-1}\\
&=n!\frac{1}{n!}\sum_{l=0}^{n}\binom{n}{l}(-1)^{n-l}(l)_{n+\alpha-k,\lambda} \nonumber \\
&=n! {n+\alpha-k \brace n}_{\lambda}.\nonumber
\end{align} \\
Hence, by \eqref{24}, \eqref{25} and \eqref{25-1}, we get 
\begin{align}
&\sum_{j=0}^{n}(-1)^{j}\binom{n}{j}(z-j)_{n+\alpha,\lambda}\label{26} \\
&=\sum_{r=0}^{n+\alpha}\binom{z-n}{r}\sum_{k=r}^{n+\alpha}\binom{n+\alpha}{k}r!{k \brace r}_{\lambda}\sum_{j=0}^{n}\binom{n}{j}(-1)^{j}(n-j)_{n+\alpha-k,\lambda}\nonumber\\
&=\sum_{r=0}^{n+\alpha}\binom{z-n}{r}\sum_{k=r}^{n+\alpha}\binom{n+\alpha}{k}r!{k \brace r}_{\lambda}n! {n+\alpha-k \brace n}_{\lambda}\nonumber \\
&=\sum_{r=0}^{n+\alpha}\binom{z-r}{r}r!n!\sum_{k=r}^{n+\alpha}\binom{n+\alpha}{k}{k \brace r}_{\lambda}{n+\alpha-k \brace n}_{\lambda}\nonumber\\
&=\sum_{r=0}^{n+\alpha}\binom{z-n}{r}r!n!\binom{n+r}{r}{n+\alpha \brace r+n}_{\lambda}\nonumber\\
&=\sum_{r=0}^{\alpha}\binom{z-n}{r}(n+r)!{n+\alpha \brace r+n}_{\lambda}.\nonumber
\end{align}
Therefore, by \eqref{26}, we obtain the following theorem. 
\begin{theorem}
For $n,\alpha,\in\mathbb{N}$, we have 
\begin{displaymath}
\sum_{j=0}^{n}(-1)^{j}\binom{n}{j}(z-j)_{n+\alpha,\lambda}=\sum_{r=0}^{\alpha}\binom{z-n}{r}(n+r)!{n+\alpha \brace r+n}_{\lambda}.
\end{displaymath}
\end{theorem}
{\it{We define $S_{n}(x,r|\lambda)$ by 
\begin{equation}
S_{n}(x,r|\lambda)=\sum_{k=0}^{n}(-1)^{k}\binom{x}{k}(k)_{r,\lambda},\quad (n,r\ge 0). \label{27}
\end{equation}}}
In particular, if $r=0$, then, by Lemma 1.1 (e), we get 
\begin{equation*}
S_{n}(x,0|\lambda)=\sum_{k=0}^{n}(-1)^{k}\binom{x}{k}=(-1)^{n}\binom{x-1}{n}. 
\end{equation*}
For $r\ge 1$, by Lemma 1.1 (d), we have 

\begin{align}
S_{n}(x,r|\lambda)&=\sum_{k=0}^{n}(-1)^{k}\binom{x}{k}(k)_{r,\lambda}=x\sum_{k=1}^{n}(-1)^{k}\binom{x}{k}\frac{k}{x}(k-\lambda)_{r-1,\lambda}\label{28} \\
&=x\sum_{k=1}^{n}(-1)^{k}\binom{x-1}{k-1}(k-\lambda)_{r-1,\lambda}\nonumber\\
&=-x\sum_{k=0}^{n-1}(-1)^{k}\binom{x-1}{k}(k+1-\lambda)_{r-1,\lambda}\nonumber\\
&=-x\sum_{k=0}^{n-1}(-1)^{k}\binom{x-1}{k}\sum_{j=0}^{r-1}\binom{r-1}{j}(k)_{j,\lambda}(1-\lambda)_{r-1-j,\lambda}\nonumber\\
&=-x\sum_{j=0}^{r-1}\binom{r-1}{j}(1-\lambda)_{r-1-j,\lambda}\sum_{k=0}^{n-1}(-1)^{k}\binom{x-1}{k}(k)_{j,\lambda}\nonumber\\
&=-x\sum_{j=0}^{r-1}\binom{r-1}{j}(1-\lambda)_{r-1-j,\lambda}S_{n-1}(x-1,j|\lambda).\nonumber
\end{align}
Therefore, by \eqref{28}, we obtain the following theorem. 
\begin{theorem}
For $r\in\mathbb{N}$, we have 
\begin{displaymath}
S_{n}(x,r|\lambda)=-x\sum_{j=0}^{r-1}\binom{r-1}{j}(1-\lambda)_{r-1-j,\lambda}S_{n-1}(x-1,j|\lambda).
\end{displaymath}
\end{theorem}
Using the definition in \eqref{27}, \eqref{4} and Lemma 1.1 (e), we note that 
\begin{align}
S_{n}(x,r|\lambda)&=\sum_{k=0}^{n}\binom{x}{k}(-1)^{k}(k)_{r,\lambda}	=\sum_{k=0}^{n}(-1)^{k}\binom{x}{k}\sum_{j=0}^{r}{r \brace j}_{\lambda}\binom{k}{j}j!\label{29} \\
&=\sum_{j=0}^{r}{r \brace j}_{\lambda}j!\sum_{k=j}^{n}(-1)^{k}\binom{x}{k}\binom{k}{j}\nonumber \\
&=\sum_{j=0}^{r}{r \brace j}_{\lambda}j!\sum_{k=j}^{n}(-1)^{k}\binom{x}{j}\binom{x-j}{k-j} \nonumber\\
&=\sum_{j=0}^{r}\binom{x}{j}{r \brace j}_{\lambda}j!(-1)^{j}\sum_{k=0}^{n-j}(-1)^{k}\binom{x-j}{k} \nonumber\\
&=\sum_{j=0}^{r}\binom{x}{j}{r \brace j}_{\lambda}j!(-1)^{j}(-1)^{n-j}\binom{x-j-1}{n-j}\nonumber\\
&=(-1)^{n}\sum_{j=0}^{r}\binom{x}{j}j!{r \brace j}_{\lambda}\binom{x-j-1}{n-j}. \nonumber
\end{align}
Therefore, by \eqref{28}, we obtain the following theorem.
\begin{theorem}
For $r\in\mathbb{N}\cup\{0\}$, we have 
\begin{equation}
S_{n}(x,r|\lambda)=(-1)^{n}\sum_{j=0}^{r}\binom{x}{j}\binom{x-j-1}{n-j}j!{r \brace j}_{\lambda}.\label{30}
\end{equation}
\end{theorem}
Replacing $x$ by $-x-1$ in \eqref{27} and using Lemma 1.1 (b), we have 
\begin{align}
S_{n}(-x-1,r|\lambda)=\sum_{k=0}^{n}(-1)^{k}\binom{-x-1}{k}(k)_{r,\lambda}=\sum_{k=0}^{n}\binom{x+k}{k}(k)_{r,\lambda}.\label{31}
\end{align}
From \eqref{6}, \eqref{30} and \eqref{31}, we note that 
\begin{align}
\sum_{k=0}^{n}\binom{x+k}{k}(k)_{r,\lambda}&=(-1)^{n}\sum_{j=0}^{r}\binom{-x-1}{j}\binom{-x-j-2}{n-j}j!{r \brace j}_{\lambda}\label{32} \\
&=(-1)^{n}\sum_{j=0}^{r}(-1)^{j}\binom{x+j}{j}(-1)^{n-j}\binom{x+n+1}{n-j}j!{r \brace j}_{\lambda}\nonumber\\
&=\sum_{j=0}^{r}\binom{x+j}{j}\binom{x+n+1}{n-j}\sum_{k=0}^{j}\binom{j}{k}(-1)^{j-k}(k)_{r,\lambda}.\nonumber
\end{align}
Therefore, by \eqref{32}, we obtain the following theorem. 
\begin{theorem}
For $r\in\mathbb{N}$, we have 
\begin{align*}
\sum_{k=0}^{n}\binom{x+k}{k}(k)_{r,\lambda}&=\sum_{j=0}^{r}\binom{x+j}{j}\binom{x+n+1}{n-j}j!{r \brace j}_{\lambda}\\
&=\sum_{j=0}^{r}\binom{x+j}{j}\binom{x+n+1}{n-j}\sum_{k=0}^{j}\binom{j}{k}(-1)^{j-k}(k)_{r,\lambda}.	\end{align*}
\end{theorem}
From Theorem 2.10, respectively with $r=0$ and $x=n$, we obtain:
\begin{align*}
&\sum_{k=0}^{n}\binom{x+k}{k}=\binom{x+n+1}{n}, \\
&\sum_{k=0}^{n}\binom{n+k}{k}(k)_{r,\lambda}=\sum_{j=0}^{r}\binom{n+j}{j}\binom{2n+1}{n-j}j!{r \brace j}_{\lambda}.
\end{align*}
By \eqref{30} and using Lemma 1.1 (b), we get 
\begin{align}
\sum_{k=0}^{n}(-1)^{k}\binom{x}{k}(k)_{r,\lambda}&=S_{n}(x,r|\lambda)=(-1)^{n}\sum_{j=0}^{r}\binom{x}{j}\binom{x-j-1}{n-j}j!{r \brace j}_{\lambda}\label{33}\\
&=\sum_{j=0}^{r}\binom{x}{j}\binom{n-x}{n-j}(-1)^{j}j!{r \brace j}_{\lambda}.\nonumber
\end{align}
\begin{corollary}
For $k\in\mathbb{N}$, we have 
\begin{equation}
\sum_{k=0}^{n}(-1)^{k}\binom{x}{k}(k)_{r,\lambda}=\sum_{j=0}^{r}(-1)^{j}\binom{n-x}{n-j}\binom{x}{j}j!{r \brace j}_{\lambda}.\label{34}	
\end{equation}
\end{corollary}
From \eqref{34} and using Lemma 1.1 (d), we note that 
\begin{align}
\sum_{k=0}^{n}(-1)^{k}\binom{x}{k}(z+yk)_{p,\lambda}&=\sum_{k=0}^{n}(-1)^{k}\binom{x}{k}\sum_{r=0}^{p}\binom{p}{r}(z)_{p-r,\lambda}(yk)_{r,\lambda}\label{35}\\
&=\sum_{k=0}^{n}(-1)^{k}\binom{x}{k}\sum_{r=0}^{p}\binom{p}{r}(z)_{p-r,\lambda}y^{r}(k)_{r,\frac{\lambda}{y}}\nonumber\\
&=\sum_{r=0}^{p}\binom{p}{r}(z)_{p-r,\lambda}y^{r}\sum_{k=0}^{n}(-1)^{k}\binom{x}{k}(k)_{r,\frac{\lambda}{y}}\nonumber\\
&=\sum_{r=0}^{p}\binom{p}{r}(z)_{p-r,\lambda}y^{r}\sum_{j=0}^{r}(-1)^{j}\binom{n-x}{n-j}\binom{x}{j}j!{r \brace j}_{\frac{\lambda}{y}}.\nonumber	
\end{align}
Therefore, by \eqref{35}, we obtain the following corollary. 
\begin{corollary}
For $n,p\ge 0$, we have 
\begin{equation}
\sum_{k=0}^{n}(-1)^{k}\binom{x}{k}(z+yk)_{p,\lambda}=\sum_{r=0}^{p}\binom{p}{r}(z)_{p-r,\lambda}y^{r}\sum_{j=0}^{r}(-1)^{j}\binom{n-x}{n-j}\binom{x}{j}j!{r \brace j}_{\frac{\lambda}{y}}. \label{36}
\end{equation}
\end{corollary} 
{\it{We define $K_{r}(x|\lambda)$ by 
\begin{equation}
K_{r}(x|\lambda)=\sum_{n=0}^{\infty}\frac{x^{n}}{n!}(n)_{r,\lambda},\quad (r\ge 0). \label{37}
\end{equation}}}
Firstly, we note that
\begin{align}
e^{x(e_{\lambda}(t)-1)}&=\sum_{j=0}^{\infty}x^{j}\frac{1}{j!}\big(e_{\lambda}(t)-1\bigl)^{j}=\sum_{j=0}^{\infty}x^{j}\sum_{r=j}^{\infty}{r \brace j}_{\lambda}\frac{t^{r}}{r!}\label{38} \\
&=\sum_{r=0}^{\infty}\sum_{j=0}^{r}x^{j}{r \brace j}_{\lambda}\frac{t^{r}}{r!}.\nonumber
\end{align}
Secondly, by the definiton of \eqref{37}, we also note that
\begin{align}
&e^{x(e_{\lambda}(t)-1)}=e^{-x}e^{xe_{\lambda}(t)} \label{39} \\
&=e^{-x}\sum_{n=0}^{\infty}\frac{x^{n}}{n!}e_{\lambda}^{n}(t)=e^{-x}\sum_{n=0}^{\infty}\frac{x^{n}}{n!}\sum_{r=0}^{\infty}(n)_{r,\lambda}\frac{t^{r}}{r!}\nonumber\\
&=\sum_{r=0}^{\infty}e^{-x}\sum_{n=0}^{\infty}\frac{x^{n}}{n!}(n)_{r,\lambda}\frac{t^{r}}{r!}=\sum_{r=0}^{\infty}e^{-x}K_{r}(x|\lambda)\frac{t^{r}}{r!}.\nonumber
\end{align}
Thus, by \eqref{39}, we get 
\begin{equation}
e^{x(e_{\lambda}(t)-1)}=\sum_{r=0}^{\infty}e^{-x}K_{r}(x|\lambda)\frac{t^{r}}{r!}.\label{40}
\end{equation}
Differentiating both sides of \eqref{40} with respect to $t$ yields the following:
\begin{equation}
xe^{x(e_{\lambda}(t)-1)}e_{\lambda}^{(1-\lambda)}(t)=\sum_{r=0}^{\infty}e^{-x}K_{r+1}(x|\lambda)\frac{t^r}{r!}. \label{41}
\end{equation}
Now, from \eqref{2}, \eqref{38}, \eqref{39}, \eqref{40} and \eqref{41}, we obtain the following result.
\begin{theorem}
For $r\ge 0$, we have 
\begin{equation}
K_{r}(x|\lambda)=e^{x}\sum_{j=0}^{r}x^{j}{r \brace j}_{\lambda}, \label{41-1}
\end{equation}
and 
\begin{displaymath}
K_{r+1}(x|\lambda)=x\sum_{j=0}^{r}\binom{r}{j}K_{j}(x|\lambda)(1-\lambda)_{r-j,\lambda}.
\end{displaymath}
\end{theorem}
Substituting $x$ by $e^{ix}=\cos x+i\sin x,\ (i=\sqrt{-1})$ in \eqref{37} and \eqref{41-1}, we obtain 
\begin{equation}
K_{r}\big(e^{ix}|\lambda\big)=\sum_{k=0}^{\infty}\frac{e^{ikx}}{k!}(k)_{r,\lambda}=e^{e^{ix}}\sum_{j=0}^{r}e^{ijx}{r \brace j}_{\lambda}. \label{42}	
\end{equation}
By \eqref{42}, we get 
\begin{align}
&\sum_{k=0}^{\infty}\frac{(k)_{r,\lambda}}{k!}\big(\cos kx+i\sin kx\big)=e^{\cos x}e^{i\sin x}\sum_{j=0}^{r}{r \brace j}_{\lambda}(\cos jx+i\sin jx)\label{43}\\
&=e^{\cos x}\sum_{j=0}^{r}{r\brace j}_{\lambda}\big(\cos(\sin x)+i\sin(\sin x)\big)\big(\cos jx+i\sin jx) \nonumber\\
&=e^{\cos x}\sum_{j=0}^{r}{r \brace j}_{\lambda}\big(\cos(\sin x)\cos jx-\sin(\sin x)\sin jx\big) \nonumber\\
&\quad +ie^{\cos x}\sum_{j=0}^{r}{r \brace j}_{\lambda}\big(\cos(\sin x)\sin jx+\sin(\sin x)\cos jx\big).\nonumber
\end{align}
By comparing the coefficients on both sides of \eqref{43}, we obtain the following theorem. 
\begin{theorem}
For $r\ge 0$, we have 
\begin{align*}
&\sum_{k=0}^{\infty}\frac{(k)_{r,\lambda}}{k!}\cos kx\\
&\quad=e^{\cos x}\bigg(\cos(\sin x)\sum_{j=0}^{r}{r \brace j}_{\lambda}\cos jx-\sin(\sin x)\sum_{j=0}^{r}{r \brace j}_{\lambda}\sin jx\bigg), 
\end{align*}
and 
\begin{align*}
&\sum_{k=0}^{\infty}\frac{(k)_{r,\lambda}}{k!}\sin kx\\
&\quad=e^{\cos x}\bigg(\cos(\sin x)\sum_{j=0}^{r}{r \brace j}_{\lambda}\sin jx+\sin(\sin x)\sum_{j=0}^{r}{r \brace j}_{\lambda}\cos jx\bigg).
\end{align*}
\end{theorem}
{\it{We define $S_{n,r}(x|\lambda)$ by
\begin{equation}
S_{n,r}(x|\lambda)=\sum_{k=0}^{n}(k)_{r,\lambda}x^{k},\quad (n,r\ge 0). \label{44}	
\end{equation}}}
The reader should compare $S_{n,r}(x|\lambda)$ here with $S_{n}(x, r|\lambda)$ in \eqref{27}. We note by using \eqref{4} and \eqref{6} that 
\begin{align*}
\sum_{k=0}^{\infty}(k)_{n,\lambda}x^{k}&=\bigg(x\frac{d}{dx}\bigg)_{n,\lambda}\bigg(\frac{1}{1-x}\bigg)=\sum_{j=0}^{n}{n \brace j}_{\lambda}x^{j}\bigg(\frac{d}{dx}\bigg)^{j}\bigg(\frac{1}{1-x}\bigg)\\
&=\sum_{j=0}^{n}{n \brace j}_{\lambda}x^{j}j!\frac{1}{(1-x)^{j+1}}=\sum_{j=0}^{n}(-1)^{j}\sum_{k=0}^{j}\binom{j}{k}(-1)^{k}(k)_{n,\lambda}\frac{x^{j}}{(1-x)^{j+1}}\\
&=\sum_{j=0}^{n}\sum_{k=0}^{j}\binom{j}{k}(-1)^{k}(j-k)_{n,\lambda}\frac{x^{j}}{(1-x)^{j+1}}, \quad (|x| <1).
\end{align*}
From \eqref{44} and using \eqref{4}, we have 
\begin{equation*}
S_{n,r}(x|\lambda)=\bigg(x\frac{d}{dx}\bigg)_{r,\lambda}\sum_{k=0}^{n}x^{k}=\bigg(x\frac{d}{dx}\bigg)_{r,\lambda}\bigg(\frac{x^{n+1}-1}{x-1}\bigg),
\end{equation*}
and 
\begin{align}
\sum_{k=0}^{n}(k)_{r,\lambda}x^{k}&=S_{n,r}(x|\lambda)=\sum_{j=0}^{r}{r \brace j}_{\lambda}x^{j}\bigg(\frac{d}{dx}\bigg)^{j}\sum_{k=0}^{n}x^{k}\label{45} \\
&=\sum_{j=0}^{r}j!{r \brace j}_{\lambda}\sum_{k=0}^{n}\binom{k}{j}x^{k}.\nonumber	
\end{align}
Therefore, by \eqref{45}, we obtain the following theorem. 
\begin{theorem}
For $n,r\ge 0$, we have 
\begin{displaymath}
\sum_{k=0}^{n}(k)_{r,\lambda}x^{k}=\sum_{j=0}^{r}j!{r \brace j}_{\lambda}\sum_{k=0}^{n}\binom{k}{j}x^{k}.
\end{displaymath}
In particular, for $x=1$, we have 
\begin{equation}
\sum_{k=0}^{n}(k)_{r,\lambda}=\sum_{j=0}^{r}j!\binom{n+1}{j+1}{r \brace j}_{\lambda}.\label{46}
\end{equation}
\end{theorem}
For $r\ge 1$ and $n\ge 0$, by Lemma 1.1 (c), Theorem 2.2 and \eqref{46}, we get 
\begin{align}
\sum_{k=0}^{n}(k)_{r,\lambda}&=\sum_{\alpha=0}^{r}\alpha!{r \brace \alpha}_{\lambda}\binom{n+1}{\alpha+1}\label{47}\\
&=\sum_{\alpha=0}^{r}\binom{n+1}{\alpha+1}\sum_{j=\alpha}^{r}(-1)^{r-j}\binom{j-1}{j-\alpha}j!{r \brace j}_{-\lambda}\nonumber\\
&=\sum_{j=0}^{r}(-1)^{r-j}{r \brace j}_{-\lambda}j!\sum_{\alpha=0}^{j}\binom{n+1}{\alpha+1}\binom{j-1}{j-\alpha}\nonumber\\
&=\sum_{j=0}^{r}(-1)^{r-j}j!{r \brace j}_{-\lambda}\sum_{\alpha=1}^{j+1}\binom{n+1}{\alpha}\binom{j-1}{j+1-\alpha}\nonumber\\
&=\sum_{j=0}^{r}(-1)^{r-j}j!{r \brace j}_{-\lambda}\sum_{\alpha=0}^{j+1}\binom{n+1}{\alpha}\binom{j-1}{j+1-\alpha}\nonumber\\
&=\sum_{j=0}^{r}(-1)^{r-j}j!{r \brace j}_{-\lambda}\binom{n+j}{j+1}. \nonumber
\end{align}
Therefore, by \eqref{47}, we obtain the following theorem. 
\begin{theorem}
For $r\ge 1$ and $n\ge 0$, we have 
\begin{equation*}
\sum_{k=0}^{n}(k)_{r,\lambda}=\sum_{j=0}^{r}(-1)^{r-j}j!{r \brace j}_{-\lambda}\binom{n+j}{j+1}.
\end{equation*}
\end{theorem}
From \eqref{46} and \eqref{47}, we get the following corollary.

\begin{corollary}
For $n\ge 0$ and $r\ge 1$, we have 
\begin{displaymath}
\sum_{j=0}^{r}j!{r \brace j}_{\lambda}\binom{n+1}{j+1}=\sum_{j=0}^{r}(-1)^{r-j}j!{r \brace j}_{-\lambda}\binom{n+1}{j+1}. 
\end{displaymath}
\end{corollary}
{\it{We define $T_{n}(x,r|\lambda)$ by
\begin{equation}
T_{n}(x,r|\lambda)=\sum_{k=0}^{n}\binom{n}{k}x^{k}(k)_{r,\lambda}.\label{50}
\end{equation}}}
We first observe that
\begin{align*}
&\Big(x\frac{d}{dx}-r \lambda \Big)_{p,\lambda}T_{n}(x,r|\lambda)=\sum_{k=0}^{n}\binom{n}{k}(k)_{r,\lambda}\Big(x\frac{d}{dx}-r \lambda \Big)_{p,\lambda}x^{k} \\
&=\sum_{k=0}^{n}\binom{n}{k}(k)_{r+p}x^{k}=T_{n}(x,r+p|\lambda).\nonumber
\end{align*}
Next, by using \eqref{4} and \eqref{6},  we have 
\begin{align}
T_{n}(x,r|\lambda)&=\Big(x\frac{d}{dx}\Big)_{r,\lambda}(1+x)^{n} \label{51}\\
&=\sum_{j=0}^{r} {r \brace j}_{\lambda}x^{j}\Big(\frac{d}{dx} \Big)^{j}(1+x)^{n} \nonumber\\
&=\sum_{j=0}^{r} {r \brace j}_{\lambda}x^{j}j!\binom{n}{j}(1+x)^{n-j} \nonumber \\
&=(1+x)^{n}\sum_{j=0}^{r}(-1)^{j}\binom{n}{j}\bigg(\frac{x}{1+x}\bigg)^{j}\sum_{k=0}^{j}(-1)^{k}\binom{j}{k}(k)_{r,\lambda}. \nonumber
\end{align}
Thus, by \eqref{51}, we get the following theorem.
\begin{theorem}
For $n\ge 0$, we have 
\begin{displaymath}
T_{n}(x,r|\lambda)=(1+x)^{n}\sum_{j=0}^{r}(-1)^{j}\binom{n}{j}\bigg(\frac{x}{1+x}\bigg)^{j}\sum_{k=0}^{j}(-1)^{k}\binom{j}{k}(k)_{r,\lambda}.
\end{displaymath}
In particular, for $x=1$, we get 
\begin{displaymath}
T_{n}(1,r|\lambda)=2^{n}\sum_{j=0}^{r}(-1)^{j}\frac{\binom{n}{j}}{2^{j}}\sum_{k=0}^{j}(-1)^{k}\binom{j}{k}(k)_{r,\lambda}.
\end{displaymath}
\end{theorem}
Now, we observe from \eqref{6} that 
\begin{align}
\big(e_{\lambda}(x)-1\big)^{n}&=\sum_{k=0}^{n}\binom{n}{k}(-1)^{n-k}e_{\lambda}^{k}(x)=\sum_{k=0}^{n}\binom{n}{k}(-1)^{n-k}\sum_{j=0}^{\infty}(k)_{j,\lambda}\frac{x^{j}}{j!}\label{55}\\
&=\sum_{j=0}^{\infty}\frac{x^{j}}{j!}\sum_{k=0}^{n}\binom{n}{k}(-1)^{n-k}(k)_{j,\lambda}\nonumber\\
&=\sum_{j=n}^{\infty}\frac{x^{j}}{j!}n!{j \brace n}_{\lambda}=\sum_{j=0}^{\infty}\frac{x^{j+n}}{(j+n)!}n!{j+n \brace n}_{\lambda}.\nonumber	
\end{align}
From \eqref{55}, we note that 
\begin{align}
\bigg(\frac{e_{\lambda}(x)-1}{x}\bigg)^{n}&=\sum_{j=0}^{\infty}\frac{n!j!}{(j+n)!}{j+n \brace n}_{\lambda}\frac{x^{j}}{j!}\label{56} \\
&=\sum_{j=0}^{\infty}\frac{{j+n \brace n}_{\lambda}}{\binom{j+n}{n}}\frac{x^{j}}{j!}.\nonumber
\end{align}
Then, by \eqref{56}, we get 
\begin{equation}
D_{x}^{n}\bigg(\frac{e_{\lambda}(x)-1}{x}\bigg)^{j}\bigg|_{x=0}=\cfrac{{n+j \brace j}_{\lambda}}{\binom{n+j}{j}},\quad (n,j \ge 0). \label{57} 
\end{equation}
By using Lemma 1.2 together with \eqref{57}, we have 
\begin{align}
&D_{x}^{n}\bigg(\frac{x}{e_{\lambda}(x)-1}\bigg)^{\alpha}\bigg|_{x=0}\label{58}\\
&=\alpha\binom{\alpha+n}{n}\sum_{j=0}^{n}(-1)^{j}\binom{n}{j}\frac{1}{\alpha+j}\bigg(\frac{x}{e_{\lambda}(x)-1}\bigg)^{\alpha+j}D_{x}^{n}\bigg(\frac{e_{\lambda}(x)-1}{x}\bigg)^{j}\bigg|_{x=0}\nonumber \\
&=\alpha\binom{\alpha+n}{n}\sum_{j=0}^{n}(-1)^{j}\binom{n}{j}\frac{1}{(\alpha+j)\binom{n+j}{n}}{n+j \brace j}_{\lambda}.\nonumber
\end{align}
On the other hand, by \eqref{8}, we get 
\begin{equation}
D_{x}^{n}\bigg(\frac{x}{e_{\lambda}(x)-1}\bigg)^{\alpha}\bigg|_{x=0}=\beta_{n,\lambda}^{(\alpha)}. \label{59}	
\end{equation}
Therefore, by \eqref{58} and \eqref{59}, we obtain the following theorem. 
\begin{theorem}
For $n\ge 0$, we have 
\begin{displaymath}
\beta_{n,\lambda}^{(\alpha)}=\alpha\binom{\alpha+n}{n}\sum_{j=0}^{n}(-1)^{j}\binom{n}{j}\frac{1}{(\alpha+j)\binom{n+j}{n}}{n+j \brace j}_{\lambda}. 
\end{displaymath}
\end{theorem}
\section{Conclusion} 
We mentioned in the Introduction that the degenerate Stirling numbers of the second kind ${n \brace k}_{\lambda}$ can be viewed as a degenerate version of the Stirling numbers of the second kind ${n \brace k}$, as the Stirling numbers of the second kind associated with the sequence ${\bf{P}}=\{(x)_{n,\lambda}\}_{n=0}^{\infty}$, and as the special case $Y=1$ of the probabilisitc degenerate Stirling numbers of the second kind associated with $Y$, $S_{2,\lambda}^{Y}(n,k)$. \par
We studied some properties, certain identities, recurrence relations and explicit expressions related to the degenerate Stirling numbers of the second kind ${n \brace k}_{\lambda}$. An explicit expression of ${n \brace k}_{\lambda}$, which is a degenerate version of Euler's formula for the Stirling numbers of the second kind, is given by
\begin{equation}
\sum_{j=0}^{k}\binom{k}{j}(-1)^{j}(j)_{n,\lambda}=(-1)^{k}k!{n \brace k}_{\lambda},\quad (n,k \ge 0). \label{60}
\end{equation}
In connection with the identity \eqref{60}, we investigated the following polynomials and power series:
\begin{align*}
&S_{n}(x, r|\lambda)=\sum_{k=0}^{n}(-1)^{k}(k)_{r, \lambda}\binom{x}{k}, \quad
K_{r}(x|\lambda)=\sum_{n=0}^{\infty}(n)_{r,\lambda}\frac{x^{n}}{n!}, \\
&S_{n,r}(x|\lambda)=\sum_{k=0}^{n}(k)_{r,\lambda}x^{k}, \quad
T_{n}(x, r|\lambda)=\sum_{k=0}^{n}\binom{n}{k}(k)_{r,\lambda}x^{k}. 
\end{align*} \par
We would like to continue to explore various degenerate versions of many special polynomials and numbers as one of our future projects.

\end{document}